\documentclass{svmult}
\usepackage{a4wide,amsfonts,amsmath,latexsym,amssymb,euscript,graphicx,mathrsfs}

\newtheorem{prop}{Proposition}[section]

\newtheorem{lemme}[prop]{Lemma}
\newtheorem{rem}[prop]{Remark}
\newtheorem{thm}[prop]{Theorem}
\newtheorem{defi}[prop]{Definition}

\newcommand{\e}{\varepsilon}
\newcommand{\numberfield}[1]{\mathbb{#1}}
\newcommand{\R}{\numberfield{R}}
\newcommand{\N}{\numberfield{N}}

\newcommand{\fin}
{ \vspace{-0.6cm}
\begin{flushright}
\mbox{$\Box$}
\end{flushright}
\noindent }

\begin{document}

\title*{{\bf A simple theory for the study of SDEs driven by
a fractional Brownian motion, in dimension one}}
\titlerunning{A simple theory for the study of SDEs driven by
a fBm, in dimension one}
\author{Ivan Nourdin}
\institute{Laboratoire de Probabilit\'es et Mod\`eles Al\'eatoires,
University Pierre et Marie Curie Paris VI,\\
Bo\^ite courrier 188, 4 Place Jussieu, 75252 Paris Cedex 5, France\\
\texttt{nourdin@ccr.jussieu.fr}
}
%
%
\maketitle

\begin{abstract}
We will
focus -- {\rm in dimension one} -- on the SDEs of the type
$dX_t=\sigma(X_t)dB_t+b(X_t)dt$ where $B$ is a fractional Brownian
motion. Our principal aim is to describe a simple
theory -- from our point of view -- allowing to study this SDE,
and this for {\rm any} $H\in (0,1)$. We will consider several definitions
of solutions and, for each of them, study conditions
under which
one has existence and/or uniqueness. Finally, we will
examine whether or not the canonical scheme associated to
our SDE converges, when the integral with respect to fBm is defined using the
Russo-Vallois symmetric integral.
\end{abstract}


\bigbreak



\section{Introduction}

The fractional Brownian motion (fBm) $B=\{B_t,\,t\ge 0\}$ of Hurst
index $H\in (0,1)$ is a centered and continuous Gaussian process
verifying $B_0=0$ a.s. and
\begin{equation}\label{fbm}
E[(B_t-B_s)^2]=|t-s|^{2H}
\end{equation}
for all $s,t\ge 0$. Observe that $B^{1/2}$ is nothing but
standard Brownian motion. Equality (\ref{fbm})
implies that the trajectories of $B$ are $(H-\e)$-H\"older
continuous, for any $\e>0$ small enough. As the fBm is selfsimilar
(of index $H$) and has stationary increments, it is used
as a model in many fields (for example, in hydrology, economics,
financial mathematics, etc.). In particular, the study of
stochastic differential equations (SDEs) driven by a
fBm is important in view of the applications. But, before
raising the question of existence and/or uniqueness for this
type of SDEs, the first difficulty is to give a meaning to the
integral with respect to a fBm. It is indeed well-known that
$B$ is not a semimartingale when $H\neq 1/2$. Thus, the It\^o or
Stratonovich calculus does not apply to this case. There are several
ways of building an integral with respect to the fBm and
of obtaining a change of variables formula. Let us point out
some of these contributions:
\begin{enumerate}
\item {\it Regularization or discretization techniques}.
Since 1993, Russo and Vallois \cite{RV93} have developed a
regularization procedure, whose philosophy is similar to the
discretization. They introduce forward (generalizing It\^o),
backward, symmetric (generalizing Stratonovich, see Definition
\ref{defsym} below) stochastic integrals and a generalized quadratic
variation. The regularization, or discretization technique, for fBm
and related processes have been performed by \cite{FP2,KZ,RV98,Z},
in the case of zero quadratic variation (corresponding to $H>1/2$).
Note also that Young integrals \cite{Y}, which are often used in this
case, coincide with the forward integral (but also with the
backward or symmetric ones, since covariation between integrand and
integrator is always zero). When the integrator has paths with
finite $p$-variation for $p>2$, 
forward and backward integrals cannot be used. 
In this case, one can use some
symmetric integrals introduced by Gradinaru {\it et al.}
in \cite{GNRV} (see $\S 2$ below). We also refer to
Errami and Russo \cite{ER} for the specific case where $H\ge 1/3$.

\item {\it Rough paths}. An other approach was taken 
by Lyons \cite{L}. His absolutely pathwise method based on 
L\'evy stochastic areas considers integrators
having $p$-variation for any $p>1$, provided one 
can construct a canonical geometric
rough path associated with the process. We refer to
the survey article of Lejay \cite{lejay} for more precise statements related to this theory.
Note however that the case where the integrator is a fBm with index $H>1/4$ has been
studied by Coutin and Qian \cite{CQ} (see also Feyel and de La Pradelle \cite{FP}).
See also Nourdin and Simon
\cite{NS2} for a link between the regularization technique and the rough paths
theory.

\item {\it Malliavin calculus}. Since fBm is a Gaussian process, it is natural
to use a Skorohod approach. Integration with respect to fBm has been
attacked by Decreusefond and \"Ust\"unel \cite{DU} for $H>1/2$ and it
has been intensively studied since (see for instance
\cite{AMN,ALN,CC}), even when the integrator is a more general
Gaussian process. We refer to Nualart's survey article 
\cite{Nual} for precise statements related to this theory.

\item {\it Wick products}. A new type of integral with zero mean defined using
Wick products was introduced by Duncan, Hu and Pasik-Duncan in \cite{DHP}, assuming
$H>1/2$. This integral turns out to coincide with the divergence operator. In \cite{bender},
Bender considers the case of arbitrary Hurst index $H\in (0,1)$ and proves an It\^o formula for
generalized functionals of $B$.
\end{enumerate}

In the sequel, we will focus -- in dimension one -- on SDEs of the type:
\begin{equation}\label{eq}
\left\{
\begin{array}{lll}
dX_t=\sigma(X_t)\,dB_t+b(X_t)dt,\,\,\,t\in [0,T]\\
X_0=x_0\in\R
\end{array}
\right.
\end{equation}
where $\sigma,b:\R\rightarrow\R$ are two continuous functions and
$H\in (0,1)$. Our principal motivation is to describe a
simple theory -- from our point of view -- allowing to study the SDE
(\ref{eq}),  for {\it any} $H\in (0,1)$.
It is linked to the regularization technique (see point 1 above).
Moreover, we emphasize that
it is already used and quoted in some research articles (see for
example \cite{BC,GN2,neuenkirch,NN,N,NS,NS2}). The aim of the
current paper is, in particular, to clarify this approach.

The paper is organized as follows. In the second part, we will
consider several definitions of solution to (\ref{eq}) and 
 for each of them we will
study under which condition one has existence
and/or uniqueness. Finally, in the third part, we will examine 
whether or not the canonical scheme associated to (\ref{eq})
converges,
when the integral with respect to fBm is defined using the
Russo-Vallois symmetric integral.

\section{Basic study of the SDE (\ref{eq})}
In the sequel, we denote by $B$ a fBm of Hurst parameter $H\in (0,1)$.
\begin{defi}
Let $X,Y$ be two real continuous processes defined on $[0,T]$. The {\rm symmetric
integral (in the sense of Russo-Vallois)} is defined by
\begin{equation}\label{defsym}
\int_0^T Y_u d^\circ X_u = 
\mathop{\rm lim\ in\ prob}\limits_{\e\rightarrow 0}
\int_0^T\frac{Y_{u+\e}+Y_u}{2}\times
\frac{X_{u+\e}-X_u}{\e}\,du,
\end{equation}
provided the limit exists and with the convention that $Y_t=Y_T$ and $X_t=X_T$ when $t>T$.
\end{defi}
\begin{rem}
{\rm
If $X,Y$ are two continuous semimartingales then $\int_0^T Y_u d^\circ X_u$ coincides
with the standard Stratonovich integral, see \cite{RV93}.
}
\end{rem}
Let us recall an important result for our study:
\begin{thm}\label{gnrv}(see \cite{GNRV}, p. 793).
The symmetric integral $\int_0^T f(B_u)d^\circ B_u$ exists for any $f:\R\rightarrow\R$
of class ${\rm C}^5$
if and only if $H\in (1/6,1)$. In this case, we have, for any 
antiderivative $F$ of $f$:
$$F(B_T)=F(0)+\int_0^T f(B_u)d^\circ B_u.$$
\end{thm}
When $H\le 1/6$, one can consider the so-called $m${\it -order
Newton-Cotes functional}:
\begin{defi}\label{def-nc}
Let $f:\R^n\rightarrow\R$ (with $n\ge 1$) be a continuous function,
$X:[0,T]\times\Omega\rightarrow\R$ and
$Y:[0,T]\times\Omega\rightarrow\R^n$ be two continuous processes and
$m\ge 1$ be an integer. The $m${\rm -order Newton-Cotes functional
of $(f,Y,X)$} is defined by
$$
\int_0^T f(Y_u) d^{{\rm NC},m} X_u = 
\mathop{\rm lim\ in\ prob}\limits_{\e\rightarrow 0}
\int_0^T \Bigl( \int_0^1 f(Y_u+\beta
(Y_{u+\e}-Y_u))\nu_m(d\beta) \Bigr)\frac{X_{u+\e}-X_u}{\e}du,
$$
provided the limit exists and with the convention that $Y_t=Y_T$ and $X_t=X_T$ when $t>T$.
Here, $\nu_1=\frac{1}{2}(\delta_0+\delta_1)$
and
\begin{equation}\label{nu}
\nu_m=\sum_{j=0}^{2(m-1)}\Bigl(\int_0^1 \prod_{k\neq j}\frac{2(m-1)u-k}{j-k}du\Bigr)\delta_{j/(2m-2)},\quad m\ge 2,
\end{equation}
$\delta_a$ being the Dirac measure at point $a$.
\end{defi}
\begin{rem}
{\rm
\begin{itemize}
\item The $1$-order Newton-Cotes functional 
$\int_0^T f(Y_u) d^{{\rm NC},1} X_u$ is nothing 
but the symmetric integral
$\int_0^T f(Y_u)d^\circ X_u$ defined by (\ref{defsym}).
On the contrary, when
$m>1$, the $m$-order Newton-Cotes functional $\int_0^T f(Y_u)
d^{{\rm NC},m} X_u$ is not {\it a priori} a``true'' integral. Indeed,
its definition could be different from $\int_0^T
\tilde{f}(\tilde{Y}_u) d^{{\rm NC},m} X_u$ even if
$f(Y)=\tilde{f}(\tilde{Y})$. This is why we call it
``{\it functional}'' instead of ``{\it integral}''.
\item The terminology ``{\it Newton-Cotes} functional'' is due to the fact that
the definition
of $\nu_m$ via~(\ref{nu}) is related to the Newton-Cotes formula of
numerical analysis. Indeed, $\nu_m$ is the unique discrete measure
carried by the numbers $j/(2m-2)$ which coincides with Lebesgue
measure on all polynomials of degree smaller than $2m-1$.
\end{itemize}
}
\end{rem}
We have the following change of variable formula.
\begin{thm}\label{gnrv-nc}(see \cite{GNRV}, p. 793).
Let $m\ge 1$ be an integer. The $m$-order Newton-Cotes functional
$\int_0^T f(B_u)d^{{\rm NC},m} B_u$ exists for any
$f:\R\rightarrow\R$ of class ${\rm C}^{4m+1}$ if and only if $H\in
(1/(4m+2),1)$. In this case, we have, for any antiderivative $F$ of $f$:
\begin{equation}\label{ito}
F(B_T)=F(0)+\int_0^T f(B_u)d^{{\rm NC},m} B_u.
\end{equation}
\end{thm}
\begin{rem}\label{rm27}
{\rm An immediate consequence of this result is that
$$\int_0^T f(B_u)d^{{\rm NC},m} B_u =\int_0^T f(B_u)d^{{\rm NC},n}
B_u=F(B_T)-F(0)$$ when $m>n$, $f$ is ${\rm C}^{4m+1}$ and $H\in
(1/(4n+2),1)$. Then, for $f$ regular enough, it is possible to
define the so-called {\it Newton-Cotes functional} $\int_0^T
f(B_u)d^{{\rm NC}} B_u$ without ambiguity by:
\begin{equation}\label{nc-gen}
\int_0^T f(B_u)d^{{\rm NC}} B_u:=
\int_0^T f(B_u)d^{{\rm NC},n} B_u\mbox{ if }H\in (1/(4n+2),1).
\end{equation}
In the sequel, we put $n_H=\inf\{n\ge 1:\,H>1/(4n+2)\}$.
An immediate consequence of (\ref{ito}) and (\ref{nc-gen}) is that, for any $H\in (0,1)$
and any $f:\R\rightarrow\R$ of class ${\rm C}^{4n_H+1}$, we have:
$$F(B_T)=F(0)+\int_0^T f(B_u)d^{{\rm NC}} B_u,$$
where $F$ is an antiderivative of $f$.
}
\end{rem}
To specify the sense of $\int_0^t \sigma(X_s)dB_s$ in
(\ref{eq}), it now seems natural to try and use the Newton-Cotes
functional. But for the time being we are only able to consider
integrands of the form $f(B)$ with $f:\R\rightarrow\R$ regular
enough, see (\ref{nc-gen}). That is why we first choose the following
definition for a possible solution to (\ref{eq}):
\begin{defi}\label{def1}
Assume that $\sigma\in {\rm C}^{4n_H+1}$ and that $b\in{\rm C}^0$. \\
{\it i)} Let $\mathfrak{C}_1$ be the class of processes
$X:[0,T]\times \Omega\rightarrow\R$ verifying that there exists
$f:\R\rightarrow\R$ belonging to ${\rm C}^{4n_H+1}$ and such that,
for every $t\in[0,T]$, $X_t=f(B_t)$ a.s.\\
{\it ii)}
A process $X:[0,T]\times \Omega\rightarrow\R$ is a solution to (\ref{eq}) if:
\begin{itemize}
\item $X\in\mathfrak{C}_1$,
\item $\forall t\in [0,T]$, $X_t=x_0+\int_0^t \sigma(X_s)d^{{\rm NC}} B_s +\int_0^t b(X_s)ds$.
\end{itemize}
\end{defi}
\begin{rem}
{\rm Note that the first point of definition {\it ii)} allows to
ensure that the integral $\int_0^t \sigma(X_s)d^{{\rm NC}}B_s$ makes
sense (compare with the adaptedness condition in the It\^o context).
}
\end{rem}
We can now state the following result.
\begin{thm}\label{thm1}
Let $\sigma\in {\rm C}^{4n_H+1}$ be a Lipschitz function, $b$ be a
continuous function and $x_0$ be a real. Then equation (\ref{eq})
admits a solution $X$ in the sense of Definition \ref{def1} if and
only if $b$ vanishes on $\mathfrak{S}(\R)$, where $\mathfrak{S}$ is
the unique solution to $\mathfrak{S}'=\sigma\circ\mathfrak{S}$ with
initial value $\mathfrak{S}(0)=x_0$. In this case, $X$ is unique and
is given by $X_t=\mathfrak{S}(B_t)$.
\end{thm}
\begin{rem}
{\rm
As a consequence of the mean value theorem, $\mathfrak{S}(\R)$ is an interval.
Moreover, it is easy to see that either $\mathfrak{S}$
is constant or $\mathfrak{S}$ is strictly monotone, and that $\inf\mathfrak{S}(\R)$
and $\sup\mathfrak{S}(\R)$ are elements of $\{\sigma=0\}\cup\{\pm\infty\}$. In particular, if $\sigma$
does not vanish, then $\mathfrak{S}(\R)=\R$ and an immediate consequence of Theorem \ref{thm1}
is that (\ref{eq}) admits a solution in
the sense of Definition \ref{def1} if and only if $b\equiv 0$.
}
\end{rem}
{\bf Proof of Theorem \ref{thm1}}. Assume that $X_t=f(B_t)$ is a solution to (\ref{eq})
in the sense of Definition \ref{def1}.
Then
\begin{equation}\label{co}
f(B_t)=x_0+\int_0^t \sigma\circ f(B_s)d^{{\rm NC}}B_s + \int_0^t b\circ f(B_s)ds=G(B_t)+\int_0^t b\circ f(B_s)ds,
\end{equation}
where $G$ is the antiderivative of $\sigma\circ f$ verifying $G(0)=x_0$.
Set $h=f-G$ and denote by $\Omega^*$ the set of $\omega\in\Omega$
such that $t\mapsto B_t(\omega)$ is differentiable at least one point
$t_0\in[0,T]$ (it is well-known that ${\rm P}(\Omega^*)=0)$. If
$h'(B_{t_0}(\omega))\neq 0$ for one $(\omega,t_0)\in\Omega\times
[0,T]$ then $h$ is strictly monotone in a neighborhood of
$B_{t_0}(\omega)$ and, for $|t-t_0|$ sufficiently small, one has
$B_t(\omega)=h^{-1}(\int_0^t b(X_s(\omega))ds)$ and, consequently,
$\omega\in\Omega^*$. Then, a.s., $h'(B_t)=0$ for all $t\in[0,T]$, so that
$h\equiv 0$. By uniqueness, one deduces $f=\mathfrak{S}$. Thus, if
(\ref{eq}) admits a solution $X$ in the sense of Definition
\ref{def1}, one necessarily has $X_t=\mathfrak{S}(B_t)$. Thanks to
(\ref{co}), one then has $b\circ\mathfrak{S}(B_t)=0$ for all $t\in
[0,T]$ a.s. and then $b$ vanishes on $\mathfrak{S}(\R)$.\fin
Consequently, when the SDE (\ref{eq}) has no drift $b$, there is a
natural solution.
But what can we do when $b\not\equiv 0$?\\
\\
Denote by $\mathscr{A}$ the set of processes $A:[0,T]\times
\Omega\rightarrow\R$ having ${\rm C}^1$-trajectories and verifying
${\rm E}\bigl(e^{\lambda \int_0^T A_s^2ds}\bigr)<\infty$ for
at least one $\lambda>1$.
\begin{lemme}\label{girsa}
Let $A\in\mathscr{A}$ and $m\in \N^*$. Then $\int_0^T
f(B_u+A_u)d^{{\rm NC},m}B_u$ exists for any $f:\R\rightarrow\R$ of
class ${\rm C}^{4m+1}$ if and only if $H>1/(4m+2)$. In this case, for any antiderivative $F$ of $f$, one has:
$$
F(B_T+A_T)=F(A_0)+\int_0^T f(B_u+A_u)d^{{\rm NC},m}B_u+\int_0^T f(B_u+A_u)A'_udu.
$$
\end{lemme}
{\bf Proof}. Set $\tilde{B}=B+A$. On the one hand, using the Girsanov
theorem in \cite{nualouk} and taking into account the assumption 
on $A$,
we have that $\tilde{B}$ is a fBm of index $H$ under 
some probability $\mathbb{Q}$
equivalent to the initial probability $\mathbb{P}$.
On the other hand, it is easy, by going back to Definition
\ref{def-nc}, to prove that $\int_0^T f(B_u+A_u)d^{{\rm NC},m}B_u$
exists if and only if $\int_0^T f(B_u+A_u)d^{{\rm NC},m}(B_u+A_u)$
does, and in this case, one has
$$
\int_0^T f(B_u+A_u)d^{{\rm NC},m}(B_u+A_u) =\int_0^T
f(B_u+A_u)d^{{\rm NC},m}B_u +\int_0^T f(B_u+A_u)A'_udu.$$ Then,
since convergence under $\mathbb{Q}$ or under $\mathbb{P}$ is
equivalent, the conclusion of Lemma \ref{girsa}
is a direct consequence of
Theorem \ref{gnrv-nc}.\fin 
Then, as previously, it is possible to
define a functional (still called {\it Newton-Cotes functional})
verifying, for any $H\in (0,1)$, for any $f:\R\rightarrow\R$ of
class ${\rm C}^{4n_H+1}$ and any process $A\in\mathscr{A}$:
$$F(B_T+A_T)=F(A_0)+\int_0^T f(B_u+A_u)d^{{\rm NC}}B_u+\int_0^T f(B_u+A_u)A'_udu,$$
where $F$ is an antiderivative of $f$.

Now, we can introduce an other definition of a  solution to (\ref{eq}):
\begin{defi}\label{def2}
Assume that $\sigma\in{\rm C}^{4n_H+1}$ and that $b\in{\rm C}^0$.\\
{\it i)} Let $\mathfrak{C}_2$ be the class of processes
$X:[0,T]\times \Omega\rightarrow\R$ such that there exist a
function $f:\R\rightarrow\R$ in ${\rm C}^{4n_H+1}$ and a process
$A\in\mathscr{A}$
such that $A_0=0$ and, for every $t\in[0,T]$, $X_t=f(B_t+A_t)$ a.s.\\
{\it ii)} A process $X:[0,T]\times \Omega\rightarrow\R$ is a solution to (\ref{eq}) if:
\begin{itemize}
\item $X\in\mathfrak{C}_2$,
\item $\forall t\in [0,T]$, $X_t=x_0+\int_0^t \sigma(X_s)d^{{\rm NC}} B_s +\int_0^t b(X_s)ds$.
\end{itemize}
\end{defi}
\begin{thm}\label{thm2}
Let $\sigma\in {\rm C}^{4n_H+1}$ be a Lipschitz function, $b$ be a continuous function
and $x_0$ be a real.
\begin{itemize}
\item If $\sigma(x_0)=0$ then (\ref{eq}) admits a solution $X$ in the sense of Definition
\ref{def2} if
and only if $b(x_0)=0$. In this case, $X$ is unique and is given by $X_t\equiv x_0$.
\item If $\sigma(x_0)\neq 0$, then (\ref{eq}) admits a solution. If moreover
$\inf_\R |\sigma|>0$ and $b\in {\rm Lip}$, this solution is unique.
\end{itemize}
\end{thm}
{\bf Proof}. Assume that $X=f(B+A)$ is a solution to (\ref{eq}) in the sense of Definition
\ref{def2}.
Then, we have
\begin{equation}\label{co2}
f(B_t+A_t)=
G(B_t+A_t)-\int_0^t \sigma(X_s)A'_sds+\int_0^t b(X_s)ds
\end{equation}
where $G$ is the antiderivative of $\sigma\circ f$ verifying $G(0)=x_0$.
As in the proof of Theorem \ref{thm1}, we obtain that
$f=\mathfrak{S}$ where $\mathfrak{S}$ is defined by
$\mathfrak{S}'=\sigma\circ\mathfrak{S}$ with initial value
$\mathfrak{S}(0)=x_0$. Thanks to (\ref{co2}), we deduce that, a.s.,
we have
$b\circ\mathfrak{S}(B_t+A_t)=\sigma\circ\mathfrak{S}(B_t+A_t)\,A'_t$
for all $t\in [0,T]$. Consequently:
\begin{itemize}
\item If $\sigma(x_0)=0$ then $\mathfrak{S}\equiv x_0$ and $b(x_0)=0$.
\item If $\sigma(x_0)\neq 0$ then $\mathfrak{S}$ is strictly monotone and
the ordinary integral equation
$$A_t=\int_0^t \frac{b\circ \mathfrak{S}}{\mathfrak{S}'}(B_s+A_s)\,ds$$
admits a maximal (in fact, global since we know already that $A$ is defined on $[0,T]$) solution by
Peano's theorem. If moreover $\inf_\R |\sigma|>0$ and $b\in {\rm Lip}$ then
$\frac{b\circ \mathfrak{S}}{\mathfrak{S}'}=\frac{b\circ \mathfrak{S}}{\sigma\circ \mathfrak{S}}\in {\rm Lip}$ and $A$ is uniquely determined.\fin
\end{itemize}
The previous theorem is not quite satisfactory because of the 
prominent role played by $x_0$. 
That is why we will finally introduce a last definition
for a solution to (\ref{eq}). We first need an analogue of Theorem
\ref{gnrv-nc} and Lemma \ref{girsa}:
\begin{thm}\label{gnrv-nc2} ({\it see \cite{these}, Chapter 4}).
Let $A$ be a process having ${\rm C}^1$-trajectories and $m\ge 1$ be
an integer. If $H>1/(2m+1)$ then the $m$-order Newton-Cotes
functional $\int_0^T f(B_u,A_u)d^{{\rm NC},m} B_u$ exists for any
$f:\R^2\rightarrow\R$ of class ${\rm C}^{2m,1}$. In this case, we
have, for any function $F:\R^2\rightarrow\R$ verifying $F'_b=f$:
$$F(B_T,A_T)=F(0,A_0)+\int_0^T f(B_u,A_u)d^{{\rm NC},m}B_u
+\int_0^T F'_a(B_u,A_u)A'_udu.$$
\end{thm}
\begin{rem}
{\rm
\begin{itemize}
\item $F'_a$ (resp. $F'_b$) means the derivative of $F$ with respect to
$a$ (resp. $b$).
\item The condition is here $H>1/(2m+1)$ and not $H>1/(4m+2)$ as
in Theorem \ref{gnrv-nc} and Lemma \ref{girsa}. Thus, for instance,
if $A\in\mathscr{A}$, if $g:\R\rightarrow\R$ is ${\rm C}^5$ and if
$h:\R^2\rightarrow\R$ is ${\rm C}^{5,1}$ then $\int_0^T
g(B_s+A_s)d^\circ B_s$ exists if (and only if) $H>1/6$ while
$\int_0^T h(B_s,A_s)d^\circ B_s$ exists {\it a priori} only when
$H>1/3$.
\item We define $m_H=\inf\{m\ge 1:\, H>1/(2m+1)\}$. As in the Remark \ref{rm27}, it is possible to
consider, for any $H\in (0,1)$ and without ambiguity, a functional
(still called {\it Newton-Cotes functional}) which verifies, for
any $f:\R^2\rightarrow\R$ of class ${\rm C}^{2m_H,1}$ and any
process $A$ having ${\rm C}^1$-trajectories:
$$F(B_T,A_T)=F(0,A_0)+\int_0^T f(B_u,A_u)d^{{\rm NC}} B_u+\int_0^T F'_a(B_u,A_u)A'_udu,$$
where $F$ is such that $F'_b=f$.
\end{itemize}
}
\end{rem}
Finally, we introduce our last definition for a solution to (\ref{eq}):
\begin{defi}\label{def3}
Assume that $\sigma\in {\rm C}^{2m_H}$ and $b\in{\rm C}^0$.\\
{\it i)} Let
$\mathfrak{C}_3$ be the class of processes $X:[0,T]\times \Omega\rightarrow\R$
verifying that there exist a function
$f:\R^2\rightarrow\R$ of class ${\rm C}^{2m_H,1}$ and a process $A:[0,T]\times \Omega\rightarrow\R$
having ${\rm C}^1$-trajectories such that $A_0=0$ and verifying,
for every $t\in[0,T]$, $X_t=f(B_t,A_t)$ a.s.\\
{\it ii)} A process $X:[0,T]\times \Omega\rightarrow\R$ is a solution to (\ref{eq}) if:
\begin{itemize}
\item $X\in\mathfrak{C}_3$,
\item $\forall t\in [0,T]$, $X_t=x_0+\int_0^t \sigma(X_s)d^{{\rm NC}} B_s +\int_0^t b(X_s)ds$.
\end{itemize}
\end{defi}
\begin{thm}\label{thm3}
Let $\sigma\in {\rm C}^{2}_b$, $b$ be a Lipschitz function
and $x_0$ be a real. Then the equation (\ref{eq}) admits a solution $X$
in the sense of Definition \ref{def3}. Moreover, if $\sigma$ is analytic, then
$X$ is the unique solution of the form $f(B,A)$ with $f$ analytic (resp.
of class ${\rm C}^1$) in the first (resp. second)
variable and $A$ a process having ${\rm C}^1$-trajectories and verifying $A_0=0$.
\end{thm}
\begin{rem}
{\rm
\begin{itemize}
\item If $H>1/3$, one can improve Theorem \ref{thm3}. Indeed, 
as shown in \cite{NS2},
uniqueness holds without any supplementary condition
on $\sigma$. Moreover, in that reference, another meaning to
(\ref{eq}) than Definition \ref{def3} is given, using the concept of
L\'evy area.
\item In \cite{NS}, one studies the problem of absolute continuity in equation (\ref{eq}),
where the solution is in the sense of Definition \ref{def3}. It
is proved that, if $\sigma(x_0)\neq 0$, then $\mathcal{L}(X_t)$ is
absolutely continuous with respect to the Lebesgue measure for all
$t\in ]0,T]$. More precisely,  the Bouleau-Hirsch
criterion is shown to
hold: if $x_t=x_0+\int_0^t b(x_s)ds$ and $t_x=\sup\{t\in
[0,T]:\,x_t\not\in {\rm Int}J\}$ where $J=\sigma^{-1}(\{0\})$ then
$\mathcal{L}(X_t)$ is absolutely continuous if and only if $t>t_x$.
\item We already said that, among the $m$-order Newton-Cotes functionals, only the first one (that is, the symmetric integral,
defined by (\ref{defsym})) is a "true" integral. For this integral, the main
results contained in this paper are summarized in the following
table (where $f$ denotes a regular enough function and $A$ a process
having ${\rm C}^1$-trajectories):
\end{itemize}
}
\end{rem}
\begin{tabular}{|c|c|c|c|c|c|}
\hline
$\begin{array}{ll}
\mbox{If we use}\\
\mbox{Definition}
\end{array}$
&
$\begin{array}{ll}
\mbox{we have}\\
\mbox{to choose}\\
H\in
\end{array}
$
&
$\begin{array}{ll}
X\mbox{ is then}\\
\mbox{of the form}
\end{array}$
&
$\begin{array}{ll}
\mbox{we have}\\
\mbox{existence if}
\end{array}$
&
$\begin{array}{ll}
\mbox{and uniqueness}\\
\mbox{if moreover}
\end{array}$
&
$\begin{array}{ll}
\mbox{See}\\
\mbox{Theorem}
\end{array}
$\\
\hline
\ref{def1}&$(1/6,1)$&$f(B)$&
$
\begin{array}{ll}
\sigma\in {\rm C}^5\cap {\rm Lip},\\
b\in{\rm C}^0\mbox{ and}\\
b_{|\mathfrak{\mathfrak{S}}(\R)}\equiv 0
\end{array}
$
&-&\ref{thm1}\\
\hline
\ref{def2}&$(1/6,1)$&$f(B+A)$&
$\begin{array}{ll}
\sigma\in {\rm C}^5\cap {\rm Lip},\\
b\in{\rm C}^0\mbox{ +}\\
i)\,\sigma(x_0)=0\\
\,\,\,\,\,b(x_0)=0\\
\mbox{or }\\
ii)\,\sigma(x_0)\neq 0
\end{array}$
&
$\begin{array}{ll}
i)\,-\\
ii)\,\inf_\R|\sigma|>0\\
\,\,\,\,\,\,\,\mbox{ and }b\in{\rm Lip}
\end{array}$
&\ref{thm2}\\
\hline
\ref{def3}&$(1/3,1)$&$f(B,A)$&
$\begin{array}{ll}
\sigma\in {\rm C}^2_b\\
\mbox{and }b\in{\rm Lip}
\end{array}$
&-&
$\begin{array}{ll}
\ref{thm3}\\
\mbox{and \cite{NS}}
\end{array}$\\
\hline
\end{tabular}
$$
\mbox{Table 1. \it{
Existence and uniqueness in SDE }}X_t=x_0+\int_0^t \sigma(X_s)d^\circ B_s + \int_0^t b(X_s)ds
$$
{\bf Proof of Theorem \ref{thm3}}. Let us remark that the classical
Doss-Sussmann \cite{D,S} method gives a natural solution $X$ of the
form $f(B,A)$. Then, in the remainder of the proof, we will
concentrate on the uniqueness. Assume that $X=f(B,A)$ is a solution
to (\ref{eq}) in the sense of Definition \ref{def3}. On the one hand, we
have
\begin{equation}\label{14}
\hskip-2cm X_t=x_0 + \int_0^t \sigma(X_s)d^{{\rm NC}}B_s
+\int_0^t b(X_s)ds
\end{equation}
$$\hskip2cm =x_0 + \int_0^t \sigma\circ f(B_s,A_s)d^{{\rm NC}}B_s
+\int_0^t b\circ f(B_s,A_s)ds.
$$
On the other hand, using the change of variables formula, we can write
\begin{equation}\label{15}
X_t=x_0+\int_0^t f'_b(B_s,A_s)d^{{\rm NC}}B_s
+\int_0^t f'_a(B_s,A_s)A'_sds.
\end{equation}
Using (\ref{14}) and (\ref{15}), we deduce that $t\mapsto\int_0^t
\varphi(B_s,A_s)d^{{\rm NC}}B_s$ has ${\rm C}^1$-trajectories where
$\varphi:=f'_b-\sigma\circ f$. As in the proof of Theorem \ref{thm2}, we
show that, a.s.,
\begin{equation}\label{16}
\forall t\in ]0,T[,\,\varphi(B_t,A_t)=0.
\end{equation}
Similarly, we can obtain that, a.s.,
$$
\forall k\in\N,\,\forall t\in ]0,T[,\,\frac{\partial^k\varphi}{\partial b^k}(B_t,A_t)=0.
$$
If $\sigma$ and $f(.,y)$ are analytic, then $\varphi(.,y)$ is analytic and
\begin{equation}\label{h=0}
\forall t\in ]0,T[,\,\forall x\in\R,\,
\varphi(x,A_t)= f'_b(x,A_t)-\sigma\circ f(x,A_t)=0.
\end{equation}
By uniqueness, we deduce
$$
\forall t\in [0,T],\,\forall x\in\R,\,f(x,A_t)=u(x,A_t),
$$
where $u$ is the unique solution to $u'_b=\sigma(u)$ with initial value $u(0,y)=y$ for any $y\in\R$.
In particular, we obtain a.s.
\begin{equation}\label{547}
\forall t\in [0,T],\,X_t=f(B_t,A_t)=u(B_t,A_t).
\end{equation}
Identity (\ref{14}) can then be rewritten as:
$$
X_t=x_0 + \int_0^t \sigma\circ u(B_s,A_s)d^{{\rm NC}}B_s
+\int_0^t b\circ u(B_s,A_s)ds,
$$
while the change of variables formula yields:
$$X_t=x_0+\int_0^t u'_b(B_s,A_s)d^{{\rm NC}}B_s
+\int_0^t u'_a(B_s,A_s)A'_sds.
$$
Since $u'_b=\sigma\circ u$, we obtain a.s.:
\begin{equation}\label{w}
\forall t\in [0,T], \,b\circ u(B_t,A_t)=u'_a(B_t,A_t)A'_t.
\end{equation}
But we have existence and uniqueness in (\ref{w}). Then the proof of Theorem is done. \fin

\section{Convergence or not of the canonical approximating schemes associated to SDE (\ref{eq})
when $d=d^\circ$} Approximating schemes for stochastic differential
equations (\ref{eq}) have already been studied only in few articles.
The first work in that direction has been proposed by Lin \cite{lin}
in 1995. When $H>1/2$, he showed that the Euler approximation of
equation (\ref{eq}) converges
uniformly in probability--but only in the easier case when 
$\sigma(X_t)$ is
replaced by $\sigma(t)$, that is, in the additive case. In 2005, I introduced in \cite{N} (see
also Talay \cite{talay}) some approximating schemes for the analogue
of (\ref{eq}) where $B$ is replaced by a H\"older continuous
function of order $\alpha$, for any $\alpha\in(0,1)$. I determined
upper error bounds and, in particular, my results apply almost
surely when the driving H\"older continuous function is a path of
the fBm $B$, for any Hurst index $H\in(0,1)$.

Results on lower error bounds are available only since very
recently: see Neuenkirch \cite{neuenkirch} for the additive case,
and Neuenkirch and Nourdin \cite{NN} (see also Gradinaru and
Nourdin \cite{GN2}) for equation (\ref{eq}). In \cite{NN}, it is
proved that the Euler scheme
$\overline{X}={\{\overline{X}^{(n)}\}}_{n\in\N^*}$ associated to
(\ref{eq}) verifies, under classical assumptions on $\sigma$ and $b$
and when $H\in(\frac{1}{2},1)$, that
\begin{equation}\label{euler-nn}
n^{2H-1}\bigl[\, \overline{X}^{(n)}_1 - X_1 \bigr]
\,{\stackrel{{\rm a.s.}}{\longrightarrow}}\,
-\frac{1}{2}\int_0^1 \sigma'(X_s)D_sX_1ds,\quad\mbox{as $n\to\infty$},
\end{equation}
where $X$ is the solution given by Theorem \ref{thm3} and $DX$ its Malliavin derivative with respect to $B$.
Still in \cite{NN}, it is proved that, for the so-called Crank-Nicholson scheme
$\widehat{X}=\{\widehat{X}^{(n)}\}_{n\in\N^*}$ associated to (\ref{eq})
with $b=0$ and defined by
\begin{equation}\label{schema}
\left\{
\begin{array}{lll}
\widehat{X}^n_0=x\\
\widehat{X}^n_{(k+1)/n}=\widehat{X}^n_{k/n}+\frac{1}{2}\bigl(\sigma(\widehat{X}^n_{k/n})
+\sigma(\widehat{X}^n_{(k+1)/n})\bigr)(B_{(k+1)/n}-B_{k/n}),\\
\hskip10cm k\in\{0,\ldots,n-1\},
\end{array}
\right.
\end{equation}
we have,
for $\sigma$ regular enough and when $H\in(\frac{1}{3},\frac{1}{2})$:
\begin{equation}\label{nn1}
\mbox{for any }\alpha<3H-1/2,\quad n^{\alpha}\bigl[ \widehat{X}^n_1 - X_1 \bigr]
\,\,{\stackrel{{\rm Prob}}{\longrightarrow}}\, \,0\mbox{ as $n\to\infty$},
\end{equation}
where $X$ is the solution given by Theorem \ref{thm1}.
Of course, this result does not give the exact rate
of convergence but only an upper bound. However, when
the diffusion coefficient $\sigma$ verifies
\begin{equation}\label{*}
\sigma(x)^2=\alpha x^2+\beta x+\gamma\mbox{ for
some $\alpha,\beta,\gamma\in\R$},
\end{equation}
the exact rate
of convergence can be derived: indeed, in this case, we have
\begin{equation}\label{nn2}
n^{3H-1/2}\bigl[ \widehat{X}^n_1 - X_1 \bigr]
\,\,{\stackrel{{\rm Law}}{\longrightarrow}}\, \,
\frac{\alpha}{12}\, \sigma(X_1)\, G,\quad\mbox{ as $n\to\infty$},
\end{equation}
with $G$ a centered Gaussian random variable independent  of $X_1$,
whose variance depends only on~$H$. Note also that, in \cite{GN2},
the exact rate of convergence associated to the schemes introduced
in \cite{N} are computed and results of the type
(\ref{nn1})-(\ref{nn2}) are obtained.

In this section, we are interested in whether
scheme (\ref{schema}) converges, according to the value of $H$ and the
expression of $\sigma$. First of all,  this problem looks
easier than computing the exact rate of
convergence, as in \cite{GN2,NN}. But, in these two papers, no
optimality is sought in the domain of validity of $H$. For
instance, in (\ref{nn1}), we impose that $H>1/3$ although it seems
more natural to only assume that $H>1/6$.

Unfortunately, we were able to find the exact barrier of convergence
for (\ref{schema}) only for particular $\sigma$, namely those which verify
(\ref{*}).
In this case, we prove in Theorem \ref{thmschema} below that the barrier
of convergence is $H=1/6$. In the other cases, it is nevertheless
possible to prove that the scheme (\ref{schema}) converge when
$H>1/3$ (see the proof of Theorem \ref{thmschema}). But the exact
barrier remains an open question.

The class (\ref{*}) is quite restricted. In particular, I must acknowledge
that Theorem \ref{thmschema} has a limited interest. However, its proof is instructive.
Moreover it contains a useful formula for $\widehat{X}^{(n)}_{k/n}$ (see Lemma \ref{lm31}), which is the core of all
the results concerning the Crank-Nicholson scheme proved in \cite{NN} (see also \cite{GN2}).

Now, we state the main result of this section:
\begin{thm}\label{thmschema}
Assume that $\sigma\in{\rm C}^1(\R)$ verifies (\ref{*}).
Then the sequence
$\{{\widehat X}^{(n)}_{1}\}$ defined by (\ref{schema}) converges in ${\rm L}^2$
if and only if $H> 1/6$.
In this case, the limit is the unique solution at time 1 to the SDE
$X_t=x_0+\int_0^t \sigma(X_s)d^\circ B_s$,
in the sense of Definition \ref{def1} and given by Theorem \ref{thm1}.
\end{thm}
\begin{rem}
{\rm
When $\sigma(x)=x$
it is easy to understand why ${\widehat X}^{(n)}_{1}$ converges in ${\rm L}^2$
if and only if $H>1/6$. Indeed, setting $\Delta_k^n=B_{(k+1)/n}-B_{k/n}$,
we have
$$
{\widehat X}^{(n)}_{1}=
x_0\prod_{k=0}^{n-1}\frac{1+\frac{1}{2}\Delta_k^n}{1-\frac{1}{2}
\Delta_k^n
}
=
x_0\,{\rm exp}\Bigl\{\sum_{k=0}^{n-1}\ln\frac{1+\frac{1}{2}\Delta_k^n}
{1-\frac{1}{2}\Delta_k^n}\Bigr\};
$$
but
$$
\ln\frac{1+\frac{1}{2}\Delta_k^n}{1-\frac{1}{2}
\Delta_k^n}=\Delta_k^n+\frac{1}{12}(\Delta_k^n)^3+\frac{1}{80}(\Delta_k^n)^5+O((\Delta_k^n)^6),
$$
and, because $\sum_{k=0}^{n-1}\Delta_k^n=B_1$ and by using Lemma \ref{lm} below, one has that
${\widehat X}^{(n)}_{1}$ converges if and only if $H>1/6$ and that,
in this case, the limit is $x_0\,{\rm exp}(B_1)$.
}
\end{rem}
As a preliminary of the proof of Theorem \ref{thmschema}, we need two lemmas:
\begin{lemme}\label{lm}
Let $m\ge 1$ be an integer.
\begin{itemize}
\item We have
$$
\sum_{k=0}^{n-1} (B_{(k+1)/n}-B_{k/n})^{2m}\mbox{ converges in }{\rm L}^2\mbox{ as }n\rightarrow\infty\mbox{ if and only if }H\ge \frac{1}{2m}.
$$
In this case, the limit is zero if $H>1/2m$ and is $(2m)!/(2^mm!)$ if $H=1/2m$.
\item We have
$$
\sum_{k=0}^{n-1} (B_{(k+1)/n}-B_{k/n})^{2m+1}\mbox{ converges in }{\rm L}^2\mbox{ as }n\rightarrow\infty\mbox{ if and only if }H>\frac{1}{4m+2}.
$$
In this case, the limit is zero.
\end{itemize}
\end{lemme}
{\bf Proof of Lemma \ref{lm}}.
The first point is an obvious consequence of the well-known convergence
$$
n^{2mH-1}\sum_{k=0}^{n-1} (B_{(k+1)/n}-B_{k/n})^{2m}
\,{\stackrel{{\rm L}^2}{\longrightarrow}}\,
(2m)!/(2^mm!),\mbox{ as }n
\rightarrow\infty.
$$
Let us then prove the second point. On the one hand, for $H>1/(4m+2)$,
we can prove directly that
$$
\sum_{k,\ell=0}^{n-1} E[(B_{(k+1)/n}-B_{k/n})^{2m+1}
(B_{(\ell+1)/n}-B_{\ell/n})^{2m+1}]
\longrightarrow 0,\mbox{ as }n\rightarrow\infty,$$
by using a Gaussian linear regression,
see for instance \cite{GRV}, Proposition 3.8. On the other hand, it is 
well known that, when $H<1/2$,
$$
n^{(2m+1)H-1/2}\sum_{k=0}^{n-1} (B_{(k+1)/n}-B_{k/n})^{2m+1}
\,{\stackrel{\mathcal{L}}{\longrightarrow}}\, {\rm N}(0,\sigma^2_{m,H}),\mbox{ as }n
\rightarrow\infty,\quad\mbox{for some $\sigma_{m,H}>0$}
$$
(use, for instance, the main result by Nualart and Peccati \cite{NP}).
We can then deduce the non-convergence when $H\le 1/(4m+2)$ as in
\cite{GNRV}, Proof of 2(c), page 796.\fin
\begin{lemme}\label{lm31}
Assume that $\sigma\in{\rm C}^5(\R)$ is bounded together with its
derivatives. Consider $\phi$ the flow associated to $\sigma$, that
is, $\phi(x,\cdot)$ is the unique solution to $y'=\sigma(y)$ with
initial value $y(0)=x$.  Then we have, for any
$\ell\in\{0,1,\ldots,n\}$:
\begin{equation}\label{final}
{\widehat X}^{(n)}_{\ell/n}=
\phi\Bigl(x_0,
B_{\ell/n}
+
\sum_{k=0}^{\ell-1}f_3({\widehat X}^{(n)}_{k/n})(\Delta_k^n)^3
+\sum_{k=0}^{\ell-1}f_4({\widehat X}^{(n)}_{k/n})(\Delta_k^n)^4
+\sum_{k=0}^{\ell-1}f_5({\widehat X}^{(n)}_{k/n})(\Delta_k^n)^5
+O(n\Delta^6(B))\Bigr).
\end{equation}
Here we set
$$
f_3=\frac{(\sigma^2)''}{24},\,
f_4=\frac{\sigma(\sigma^2)'''}{48}
\mbox{ and }
f_5=
\frac{\sigma'^4}{80}+\frac{\sigma^2\sigma'\sigma'''}{15}+\frac{3\sigma\sigma'^2\sigma''}{40}
+\frac{\sigma^2\sigma''^2}{20}+\frac{\sigma^3\sigma^{(4)}}{80},
$$
$$
\Delta_k^n=B_{(k+1)/n}-B_{k/n},\mbox{ when }n\in\N\*\mbox{ and }k\in\{0,1,\ldots,n-1\}
$$
and
$$
\Delta^{p}(B) = \max_{k=0, \ldots, n-1}
\left|\left(\Delta_k^{n}\right)^p\right|,\mbox{ when $p \in \N^*$.}
$$
\end{lemme}
{\bf Proof of Lemma \ref{lm31}}.
Assume, for an instant, that $\sigma$ does not vanish.
In this case, $\phi(x,\cdot)$ is a bijection from $\R$ to himself for any $x$ and we can consider $\varphi(x,\cdot)$
such that
\begin{equation}\label{phi}
\forall x,t\in\R:\,\varphi(x,\phi(x,t))=t\mbox{ and }
\phi(x,\varphi(x,t))=t.
\end{equation}
On the one hand, thanks to (\ref{phi}), it is a little long but easy to compute
that
$$
\begin{array}{llllll}
\varphi(x,x)&=&0,\\
\varphi'_t(x,x)&=&1/\sigma(x),\\
\varphi''_{tt}(x,x)&=&[-\sigma'/\sigma^2](x),\\
\varphi^{(3)}_{ttt}(x,x)&=&[(2\sigma'^2-\sigma\sigma'')/\sigma^3](x),\\
\varphi^{(4)}_{tttt}(x,x)&=&[(-6\sigma'^3+6\sigma\sigma'\sigma''
-\sigma^2\sigma''')/\sigma^4](x)\\
\varphi^{(5)}_{ttttt}(x,x)&=&[(24\sigma'^4-36\sigma\sigma'^2\sigma''+8\sigma^2\sigma'\sigma'''
+6\sigma^2\sigma''^2-\sigma^3\sigma^{(4)})/\sigma^5](x).
\end{array}
$$
Then, for $u$ sufficiently small, we have
$$
\begin{array}{lll}
\varphi(x,x+u)&=\frac{1}{\sigma}(x)u-\frac{\sigma'}{2\sigma^2}(x)u^2+
\frac{2\sigma'^2-\sigma\sigma''}{6\sigma^3}(x)u^3+\frac{-6\sigma'^3+6\sigma\sigma'\sigma''-
\sigma^2\sigma'''}{24\sigma^4}(x)u^4\\
&\hskip3.1cm{}+\frac{24\sigma'^4-36\sigma\sigma'^2\sigma''+8\sigma^2\sigma'\sigma'''+6\sigma^2\sigma''^2
-\sigma^3\sigma^{(4)}}{\sigma^5}(x)u^5+O(u^6).
\end{array}
$$
On the other hand, using (\ref{schema}) and some basic Taylor expansions,
one has for $k\in\{0,1,...,n-1\}$:
$$
\begin{array}{lll}
{\widehat X}^{(n)}_{(k+1)/n}&=&
{\widehat X}^{(n)}_{k/n}
+\sigma({\widehat X}^{(n)}_{k/n})
\Delta_k^n
+\frac{\sigma\sigma'}{2}({\widehat X}^{(n)}_{k/n})(\Delta_k^n)^2
+\frac{\sigma\sigma'^2+\sigma^2\sigma''}{4}({\widehat X}^{(n)}_{k/n})(\Delta_k^n)^3\\
&&{}+\bigl(
\frac{\sigma\sigma'^3}{8}+\frac{3\sigma^2\sigma'\sigma''}{8}+\frac{\sigma^3\sigma'''}{12}
\bigr)({\widehat X}^{(n)}_{k/n})(\Delta_k^n)^4\\
&&{}+\bigl(
\frac{\sigma\sigma'^4}{16}+\frac{3\sigma^2\sigma'^2\sigma''}{8}+\frac{\sigma^3\sigma'\sigma'''}{6}
+\frac{\sigma^3\sigma''^2}{8}+\frac{\sigma^4\sigma^{(4)}}{48}
\bigr)({\widehat X}^{(n)}_{k/n})(\Delta_k^n)^5\\
&&{}+O(\Delta^6(B)).
\end{array}
$$
Then, we have
$$
\begin{array}{lll}
\varphi(
{\widehat X}^{(n)}_{k/n}
,
{\widehat X}^{(n)}_{(k+1)/n}
)&=&
\varphi({\widehat X}^{(n)}_{k/n},{\widehat X}^{(n)}_{k/n}
+[{\widehat X}^{(n)}_{(k+1)/n}-{\widehat X}^{(n)}_{k/n}])\\
&=&
\Delta_k^n
+\frac{\sigma'^2+\sigma\sigma''}{12}({\widehat X}^{(n)}_{k/n})(\Delta_k^n)^3
+\bigl(\frac{\sigma\sigma'\sigma''}{8}+
\frac{\sigma^2\sigma'''}{24}
\bigr)({\widehat X}^{(n)}_{k/n})(\Delta_k^n)^4 \\
&&{}+\bigl(
\frac{\sigma'^4}{80}+\frac{\sigma^2\sigma'\sigma'''}{15}+\frac{3\sigma\sigma'^2\sigma''}{40}
+\frac{\sigma^2\sigma''^2}{20}+\frac{\sigma^3\sigma^{(4)}}{80}
\bigr)({\widehat X}^{(n)}_{k/n})(\Delta_k^n)^5
+O(\Delta^6(B))\\
&=&
\Delta_k^n
+f_3
({\widehat X}^{(n)}_{k/n})(\Delta_k^n)^3
+f_4
({\widehat X}^{(n)}_{k/n})(\Delta_k^n)^4
+f_5({\widehat X}^{(n)}_{k/n})(\Delta_k^n)^5
+O(\Delta^6(B)).
\end{array}
$$
We deduce, using (\ref{phi}):
$$
{\widehat X}^{(n)}_{(k+1)/n}=
\phi\bigl({\widehat X}^{(n)}_{k/n},
\Delta_k^n
+f_3
({\widehat X}^{(n)}_{k/n})(\Delta_k^n)^3
+f_4
({\widehat X}^{(n)}_{k/n})(\Delta_k^n)^4
+f_5({\widehat X}^{(n)}_{k/n})(\Delta_k^n)^5
+O(\Delta^6(B))
\bigr).
$$
Finally, by using the semi-group property verified by $\phi$, namely
$$
\forall x,s,t\in\R:\,\phi(\phi(x,t),s)=\phi(x,t+s).
$$
we easily deduce (\ref{final}).

In fact, we assumed that $\sigma$ does not vanish only for having
the possibility to introduce $\varphi$. But (\ref{final}) is an
algebraic formula then it is also valid for general $\sigma$, as
soon as it is bounded together with its derivatives. \fin {\bf Proof
of Theorem \ref{thmschema}}. Assume that $\sigma$ verifies
(\ref{*}). Although $\sigma$ is not bounded in general, it is easy
to verify that we still have $O(n\Delta^6(B))$ as remainder in
(\ref{final}). Moreover, simple but tedious computations show that
we can simplify in (\ref{final}) to obtain
$$
\begin{array}{lll}
{\widehat X}^{(n)}_{1}=
\phi\bigl(x_0,
B_{1}+
\frac{\alpha}{12}
\sum_{k=0}^{n-1}
(\Delta_k^n)^3
+\frac{\alpha^2}{80}\sum_{k=0}^{n-1}(\Delta_k^n)^5
+O(n\Delta^6(B))\bigr).
\end{array}
$$
Thus, as a conclusion of Lemma \ref{lm}, we obtain easily that
${\widehat X}^{(n)}_{1}$ converges to $\phi(x_0,B_1)$ if and only if
$H>1/6$. \fin
 \noindent
\small {\bf Acknowledgement}. I am indebted to the anonymous referee
for the careful
 reading of the original manuscript and for a number of suggestions.

\end{document}